\documentclass{amsart} %\usepackage{amsfonts}\usepackage{amssymb}   

\usepackage[ref]{overcit} %modified overcite.sty to put paren

\newcommand{\lbl}[1]{{\tt \small [#1]}\label{#1}}
\renewcommand{\lbl}[1]{\label{#1}}

\newcommand{\omitted}[1]{{\em ... omitted text ...}}
 
\newcommand{\mmm}{\hat{\mu}_n}
\newcommand{\nnn}{\hat{\nu}_n}

\newcommand{\rf}[1]{(\ref{#1})}

\newcommand{\calP}{{\mathcal P}}
\newcommand{\calK}{{\mathcal K}}

\newcommand{\sR}{{\mathbb R}}

\newcommand{\II}{{\mathbb I}}
\newcommand{\LL}{{\mathbb L}}
\newcommand{\KK}{{\mathbb K}}

\newcommand{\xx}{{\mathbf x}}
\newcommand{\yy}{{\mathbf y}}

\newcommand{\eps}{\varepsilon}

       \newtheorem{Theorem}{Theorem}%[section]
       \newtheorem{Proposition}{Proposition}
       \newtheorem{Corollary}{Corollary} 
        
       \newtheorem{Lemma}{Lemma} 
       \newtheorem{Remark}{Remark}
\newtheorem{Claim}{Claim} 
\theoremstyle{definition}
\newtheorem{Definition}{Definition}  
 
\newtheorem{Assumption}{Assumption}

\newenvironment{proofof}[1]{\noindent {\bf #1.\/}}{\qed\vskip 0.1in}
\author{
W{\l}odzimierz  Bryc
}
\thanks{Department of Mathematics,
University of Cincinnati,
P.O. Box 210025,
Cincinnati, OH 45221--0025. Email: 
Wlodzimierz.Bryc@UC.edu}
\address{Department of Mathematics \\
University of Cincinnati\\
P.O. Box 210025\\
Cincinnati, OH 45221--0025\\
Wlodzimierz.Bryc@UC.edu}

\keywords{large deviations, symmetric interaction, non-convex rate function}
\subjclass[2000]{ Primary: 60F10 Secondary: 60K35}
\date{June 21, 2003} %second date fix: May 15, 2002% first date fixed at: March 17, 2002

\title[Large deviations under  symmetric interactions]
{Large deviations of  empirical measures under symmetric interaction}
\begin{document}
\begin{abstract}
We prove the large deviation principle for the joint empirical measure of pairs of random
variables which are coupled by a ``totally symmetric" interaction. The
rate
function is given by an explicit bilinear expression, 
which is finite only on product measures and hence is
non-convex. 

\end{abstract}
\maketitle
\tolerance=7000
 
\section{Introduction}
\subsection{}%Motivation, History, References}

Large deviations of empirical measures have been widely studied in the literature since the
celebrated Sanov's theorem, %\cite{Sanov-LD}, \cite{Sanov-LD-tr}, 
which gives  the large
deviations principle in the scale of $n$ of the empirical measures of i.~i.~d. random variables
with the  relative entropy $H(\mu|\nu)=\int \log \frac{d\mu}{d\nu}d\mu$ as the rate function. 
%The
%relative entropy emerges as a component of a rate function in many other situations, including
%particle systems with interactions of finite range, see for example \cite{MR2001m:60060}.
%%\cite{BenArous-Guionnet-97a},
%\cite{MR2001m:60060},
%\cite {Moral-Guionnet-98},
%\cite{DelMoral-Guionnet-01}.
%%\cite{BenArous-Guionnet-95},
%%\cite{BenArous-Guionnet-95e},
%%\cite{BenArous-Guionnet-98}
%\cite{Timo93a},\cite{Timo93b},\cite{Timo94}, \cite{Timo95b}, \cite{Timo95a},
% \cite{Timo98b}, \cite{Timo98a}.
%\cite{BenArous-Dembo-Guionnet-01}?
%
Another entropy, Voiculescu's non-commutative entropy $\Sigma(\mu)=\iint \log |x-y|
\mu(dx)\mu(dy)$, arises in the study of fluctuations of eigenvalues of
random  matrices, see Hiai \& Petz\cite{Hiai-Petz-00} 
and the references therein. 
%summarizes new results, including
% \cite{BenArous-Guionnet-97}
%\cite{BenArous-Zeitouni}
% \cite{Johansson-98}. 
%Analogies
% between the ``non-commutative" and commutative entropy have been
%studied for example in \cite{Biane-01},  \cite{Hiai-Petz98}, \cite{Hiai-Petz-00}, \cite{Voiculescu-99}. 
Chan \cite{Chan} interprets empirical measures of eigenvalues of random matrices 
 as a system of interacting diffusions with singular interactions.

\subsection{}%This paper}
In this paper we   study empirical measures which can be thought of as a decoupled version of the empirical
measures generated by random matrices.  We are interested in empirical
measures on $\sR^2$ generated
by pairs of random variables that are tied together by a totally symmetric, and
hence non-local, interaction, see formula \rf{f} for the (unnormalized) joint density. 
Under certain assumptions, we prove that the large deviation principle in the scale $n^2$ holds
for the joint empirical measures, and the rate function is  non-convex.  As a corollary,
we derive a large deviations principle for the univariate average 
empirical measures with a rate function that superficially
resembles
the rate function of random matrices, see Corollary \ref{C2}; 
an interesting feature here is the emergence of concave rate functions, see Remark \ref{RC2}. 
(Eigenvalues of random matrices are exchangeable and the large deviation rate function for their
 empirical measures is convex;
 infinite exchangeable sequences often lead to
 non-convex rate functions, see Dinwoodie \& Zabell \cite{Dinwoodie-Zabell} and \cite[Example 3]{Bryc90a}.) 
%(The assumptions of Corollary \ref{C2} exclude the integral kernel which arises in the study of
%eigenvalues of random matrices.) 

\subsection{}%The setup}
Let $g:\sR^2\to\sR$ be a continuous function which satisfies the following conditions. 
\begin{Assumption}\lbl{B0} $g(x,y)\geq 0$ for all $x,y\in\sR$.
\end{Assumption}
\begin{Assumption}%[Forgotten assumption]
\lbl{p-int} For every $0<\alpha\leq 1$,
$M_\alpha:=\iint g^\alpha(x,y) dx dy<\infty$.
\end{Assumption}
\begin{Assumption}\lbl{B1} $g(x,y)$ is bounded, $g(x,y)\leq e^C$.
\end{Assumption}
In the following statements we use the convention that $-\log 0=\infty$. 
\begin{Assumption}\lbl{B2} The function $k(x,y):=-\log g(x,y)$ has compact level sets:
for every $a>0$ the set $\{(x,y): g(x,y)\geq e^{-a}\}\subset\sR^2$ is compact.
\end{Assumption}

The purpose of the next assumption is to allow singular interactions, where
$g(x,x)=0$; this assumption is 
automatically satisfied with 
$\beta=0$ if $g(x,y)>0$ for all $x,y$.
\begin{Assumption}\lbl{B3} There is a $\beta\geq 0$ such that 
$(x,y)\mapsto\beta\log|x-y|-\log g(x,y)$ extends from $\{(x,y): x\ne y\}$ to the  continuous function on $\sR^2$.
\end{Assumption}

Examples of functions that satisfy these assumptions are: the Gaussian kernel 
$$g(x,y)=e^{-x^2-y^2+2\theta xy}$$ for $|\theta|<1$, see the proof of
Proposition \ref{IIG2};
and a singular kernel $$g(x,y)=|x-y|^\beta e^{-x^2-y^2}$$  for
 $\beta\geq 0$, see the proof of Proposition \ref{LOG}.

%\subsubsection{Probability measure}
Define
\begin{equation}\lbl{f}
f(x_1,\dots,x_n,y_1,\dots,y_n)=\prod_{i,j=1}^n g(x_i,y_j).
\end{equation}
Clearly, $f$ depends on $n$; we will suppress this dependence in our notation and we will
further write
$f(\xx,\yy)$ as a convenient shorthand for $f(x_1,\dots,x_n,y_1,\dots,y_n)$. 

Assumptions \ref{B0}, \ref{p-int}, and \ref{B1} imply that $f$ is integrable. Indeed, 
since $g(x,y)\leq e^{C}$, 
$$Z_n:=\int_{\sR^{2n}} f(x_1,\dots,x_n,y_1,\dots,y_n)\,d x_1\dots dx_n dy_1\dots dy_n
$$
$$\leq \int_{\sR^{2n}} \prod_{i=1}^n \left(g(x_i,y_i)e^{C(n-1)}\right)\,d x_1\dots dx_n dy_1\dots dy_n
= e^{C{(n^2-n)}}M_1^n <\infty.$$
 
We are interested in  
joint empirical measures 
\begin{equation}\lbl{mu_n}
\mmm=\frac{1}{n^2}\sum_{i,j=1}^n\delta_{x_i,y_j},
\end{equation}
considered as random variables with values in the Polish space of probability measures
$\calP(\sR^2)$ (equipped  with the topology of weak convergence),   with the distribution induced
on
$\calP(\sR^2)$
 by the probability measure $\Pr=\Pr_n\in\calP(\sR^{2n})$ defined by
\begin{equation}\lbl{Pr}
\Pr(d\xx,d\yy):=\frac{1}{Z_n}f(\xx,\yy)\,d\xx
d\yy.
\end{equation}

%\subsection{The main theorem}
 
\begin{Theorem} \lbl{T1}
%If $k(x,y)$ satisfies assumptions \ref{A0}, \ref{A1},
If $g(x,y)$ satisfies Assumptions \ref{B0}, \ref{p-int}, \ref{B1}, \ref{B2}, and \ref{B3}, 
 then the joint empirical measures
$\{\mmm\}$ satisfy the large deviation principle in the scale $n^2$ with the rate function
$\II:\calP(\sR^2)\to [0,\infty]$ given by
\begin{equation}\lbl{I}
\II(\mu)=\left\{\begin{array}{ll}
\iint k(x,y)\nu_1(dx)\nu_2(dy) -I_0 &
\mbox{if $\mu=\nu_1\otimes\nu_2$  is a product}\\
& \mbox{ measure and $k$ is $\mu$-integrable;} \\
& \\
\infty &\mbox{otherwise,}
\end{array}\right.
\end{equation}
 where $k(x,y)=-\log g(x,y)$ and 
$I_0=\inf_{x,y\in\sR} k(x,y)$.
\end{Theorem}

\begin{Definition}[{\cite[Chapter 3]{BCR}}] We say that $k:\sR^2\to\sR$ is a negative definite kernel 
if $k(x,y)=k(y,x)$ and 
\begin{equation}\lbl{neg-def}
\sum k(x_i,x_j)c_ic_j\leq 0
\end{equation} for 
all $x_i,c_i\in\sR$ such that $\sum c_i=0$.
\end{Definition}
Condition \rf{neg-def} is satisfied for
 $k(x,y)=V(x)+W(y)-\kappa(x,y)$, where $\kappa(x,y)$ is 
positive-definite.

Consider the average empirical measures
$$\hat{\sigma}_n:=\frac1{2n}\sum_{i=1}^n (\delta_{x_i}+\delta_{y_i}).$$

\begin{Corollary}\lbl{C2}  Suppose that
 the assumptions of Theorem \ref{T1} hold true, and in addition $k(x,y)$ is continuous   and
negative-definite.    Then the average
empirical measures $\{\hat{\sigma}_n\}$
satisfy the large deviation principle  in the scale $n^2$ with the rate function
$$\II(\nu)=\iint k(x,y) \nu(dx)\nu(dy)-I_0,$$
and $I_0=\inf_x k(x,x)$.
\end{Corollary}

\begin{proof}This follows from the contraction principle. 
The mapping 
$\mu()\mapsto \frac12\int \mu(\cdot, dy)+ \frac12\int \mu(dx, \cdot)$ is continuous in the weak topology. The
rate function is
$\II(\nu)=\inf \{\II(\nu_1\otimes\nu_2): \nu=\frac12\nu_1+\frac12\nu_2\}$. 

Write $\KK(\mu)=\iint k(x,y)\mu(dx,dy)$. 
If $\nu=\frac12\nu_1+\frac12\nu_2$ then 
%by  the standard approximation argument
%\cite[Chapter 7, Proposition 1.2]{BCR} 
\cite[Theorem 3]{Ressel}
implies that 
$\KK(\nu_1\otimes\nu_2)\geq \KK(\nu\otimes\nu)$. %, cf. \cite[Chapter 7, Proposition 1.2]{BCR}.
Thus  
$\II(\nu)=\KK(\nu\otimes\nu)-I_0$. 

Another form of the cited inequality is that for any two probability measures $\nu_1,\nu_2$ we have
\begin{equation}\lbl{2K}2\KK(\nu_1\otimes\nu_2)\geq \KK(\nu_1\otimes\nu_1)+\KK(\nu_2\otimes\nu_2).
\end{equation}
In particular,  $2k(x,y)\geq k(x,x)+k(y,y)$, which implies that
$I_0=\inf_{x,y} k(x,y)=\inf_x k(x,x)$. 
%\mybox{ZW 61 (1982) 223--235} \mybox{ZW 57 (1981) 193--201}
\end{proof} 

\begin{Remark}\lbl{RC1} Inequality \rf{2K}  implies that
 the rate function satisfies %an unusual inequality:
$$\II\left(\frac12\nu_1+\frac12\nu_2\right)\geq
\frac12\II(\nu_1)+\frac12\II(\nu_2).$$
\end{Remark}

\section{Applications}

\subsection{}
%This section illustrates that \rf{f} may lead to  natural rate functions.
Let $g(x,y)=e^{-x^2-y^2+2\theta xy}$.
Then
$$f(\xx,\yy)=\exp(-n\sum_{i=1}^nx_i^2-n\sum_{j=1}^ny_j^2+2\theta\sum_{i,j=1}^nx_iy_j)$$
and $k(x,y)=x^2+y^2-2\theta xy$. 

Denote by $m_r(\nu)=\int x^r\nu(dx)$  the $r$-th moment of a measure $\nu$.

\begin{Proposition}\lbl{IIG2}
\begin{itemize}
\item[(i)]
If $|\theta|<1$  then the empirical measures
$$\nnn:=\frac1n\sum_{i=1}^n \delta_{x_j}$$
satisfy the large deviation principle in the scale $n^2$ with the rate function 
$$\II(\nu)=
\left(m_2(\nu)-\theta^2 m_1^2(\nu)\right).$$
\item[(ii)]
 If $0\leq \theta<1$  then the average 
empirical measures
$$\hat{\sigma}_n:=\frac1{2n}\sum_{i=1}^n (\delta_{x_i}+\delta_{y_i}).$$ 
satisfy the large deviation principle  in the scale $n^2$ with the rate function
$$\II(\nu)=2\left(m_2(\nu)-\theta m_1^2(\nu)\right).$$
\end{itemize}
(In the formulas above, use $\II(\nu)=\infty$ if $m_2(\nu)=\infty$.)
\end{Proposition} 
\begin{Remark}The marginal density  relevant in {Proposition} \ref{IIG2}(i) is
$$f_1(\xx)=C(n,\theta) \exp\left({-n^2\left(\frac1n\sum_{i=1}^n x_i^2-\theta^2(\frac1n\sum_{i=1}^n
x_i)^2\right)}\right).$$
%compare \cite[(2.1)]{BenArous-Dembo-Guionnet-01}.
\end{Remark}
\begin{Remark}\lbl{RC2} Both rate functions in {Proposition} \ref{IIG2} are concave.
\end{Remark}
%\mybox{What if $\theta<0$?}
\begin{proof}(i)
It is easy to see that the assumptions of Theorem \ref{T1} are satisfied.
Indeed, $k(x,y)=(x-\theta y)^2+(1-\theta^2)y^2$ is continuous, bounded from below. Furthermore
$$\{k(x,y)\leq a^2\}\subset \{|y|\leq |a|/(1-\theta^2)\}\cap\{ |x|\leq |a|/(1-\theta^2)\}$$ so $k(x,y)$ has
compact level sets. Finally, for $\alpha>0$ by a change of variables we see that  
$\iint e^{-\alpha k(x,y)}dx dy =\frac{1}{\alpha} \iint e^{-k(x,y)}dx dy<\infty$ so Assumption \ref{p-int} is satisfied, too. 

The result follows by the contraction principle: taking a marginal of a measure in $\calP(\sR^2)$ is a
continuous mapping. The rate function is
$\inf \{\II(\mu): \nu(A)=\mu(A\times \sR)\}$. But since $\II$ is infinite on non-product measures,
this is the same as $\inf_{\nu_2}\{\iint k(x,y)\nu(dx)\nu_2(dy)-I_0\}$.
 Since $I_0=0$ here, it remains to notice that 
$\inf_{\nu_2}\{\iint k(x,y)\nu(dx)\nu_2(dy)\}=\inf_{y}\{\int k(x,y)\nu(dx)\}=
\inf_y \{m_2(\nu)+y^2-2\theta y m_1(\nu)\}=m_2(\nu)-\theta^2m_1^2(\nu)$.

(ii) This follows  from Corollary \ref{C2}:
 if $\theta\geq 0$ then $2\theta xy$ is positive-definite. Thus 
$k(x,y)=x^2+y^2-2\theta xy$ is a negative definite kernel.
%satisfies condition \rf{neg-def}
%$$\sum k(x_i,x_j)\la_i\la_j=-2\theta\left(\sum_j\la_jx_j\right)^2.$$
\end{proof}

\subsection{}%Singular interactions}
Next, we consider a model which can be interpreted as a ``decoupled" version of a model 
studied in relation to eigenvalue fluctuations of random matrices, where one
encounters $x_j$ instead of our $y_j$, compare 
\cite[Section 5]{BenArous-Guionnet-97}, \cite[formula (1.9)]{Johansson-98}. 
We consider here a slightly more general situation when
 $$g(x,y)=|x-y|^\beta e^{-V(x)-W(y)}.$$ Then
$$f(\xx,\yy)=\prod_{i,j=1}^n|x_i-y_j|^\beta \prod_{i=1}^ne^{-nV(x_i)}\prod_{j=1}^ne^{-nW(y_j)},$$ 
and 
 $k(x,y)=V(x)+W(y)-\beta\log|x-y|$.
We assume that functions $V(x), W(y)$ are continuous, $\beta\geq 0$, and that  
\begin{equation}\lbl{VW}
\lim_{|x|\to\infty} \frac{V(|x|)}{\log\sqrt{1+x^2}}=
\lim_{|y|\to\infty} \frac{W(|y|)}{\log\sqrt{1+y^2}}=\infty.
\end{equation}

\begin{Proposition}\lbl{LOG} The bivariate empirical measures 
$\mmm$ defined by \rf{mu_n} satisfy the large deviation principle in the scale $n^2$
with the rate function $\II$ given by \rf{I}.

In particular, if $V(u)=W(u)=u^2$, then the rate function is
$$\II(\nu_1\otimes\nu_2)=m_2(\nu_1)+m_2(\nu_2)-\beta \iint \log|x-y|\nu_1(dx)\nu_2(dy)+1/2\ \beta(\log \beta - 1).$$

\end{Proposition}
\begin{proof}
We verify that the hypotheses of Theorem \ref{T1} are satisfied.
Assumption \ref{B0} holds trivially. Assumption \ref{B3} holds trivially since $V(x)+W(y)$ is
continuous.

To verify Assumption \ref{B1} notice that
\begin{equation}\lbl{k>V+W}
k(x,y)\geq V(x)+W(y)-\beta\log\sqrt{1+x^2}-\beta\log\sqrt{1+y^2}.
\end{equation}

Since  $V(x)-\beta\log\sqrt{1+x^2}$ is a continuous function which by \rf{VW} tends to infinity
as $x\to\pm \infty$,
 it is bounded from below, $V(x)-\beta\log\sqrt{1+x^2}\geq -c$ for
some $c$. Similarly, $W(y)-\beta\log\sqrt{1+y^2}\geq -c$.

We now verify Assumption \ref{B2}. The set  $K_a:=\{k(x,y)\leq a\}$ is 
closed since $k$ is lower semicontinuous. Furthermore,
\rf{k>V+W} implies that $K_a$ is contained in a level set of the continuous function 
$V(x)+W(y)-\beta\log\sqrt{1+x^2}-\beta\log\sqrt{1+y^2}$. The latter  set is bounded
since $V(x)-\beta\log\sqrt{1+x^2}>a+c$ for all large enough $|x|$ and similarly  
$W(y)-\beta\log\sqrt{1+y^2}>a+c$ for all large enough $|y|$. 

To verify Assumption \ref{p-int} we use
inequality \rf{k>V+W}  again. It implies
 $$\iint g(x,y)^\alpha dx dy\leq \int e^{-\alpha (V(x)-\beta\log\sqrt{1+x^2})}dx\int e^{-\alpha (W(y)-\beta\log\sqrt{1+y^2})}dy.$$
By assumption \rf{VW}, there is $N>0$ such that
 for $|x|>N$ we have $V(x)>(\beta+2/\alpha) \log\sqrt{1+x^2}$. 
By the previous argument the integrand is bounded; thus
$\int e^{-\alpha (V(x)-\beta\log\sqrt{1+x^2})}dx
\leq \int_{-N}^{N} e^{-\alpha (V(x)-\beta\log\sqrt{1+x^2})}dx+
\int_{|x|>N}e^{-\alpha (2/\alpha\log\sqrt{1+x^2})}dx \leq 2N e^{\alpha c}+\int_{|x|>N}\frac{1}{1+x^2}dx<\infty$.

Therefore, by Theorem \ref{T1} the  empirical measures $\mmm$ satisfy the large deviation principle  with the rate function 
$\II(\nu_1\otimes\nu_2)= \int V(x) \nu_1(dx)+\int W(y) \nu_2(dy) -\beta \iint \log |x-y|
\nu_1(dx)\nu_2(dy)-I_0$.
If $W(u)=V(u)=u^2$ then $I_0=\inf_{x,y}\{x^2+y^2-\beta\log|x-y|\}=\beta/2(1-\log \beta)$ by calculus.
\end{proof}

\section{Auxiliary results and proof of Theorem {\protect\ref{T1}}}
The proof relies on Varadhan's functional method, see 
\cite[Theorem T.1.3]{Bryc90a}, \cite[Theorem 4.4.10]{Dembo-Zeitouni}.
It consists of two steps: verification that the Varadhan  functional
$$\Phi\mapsto \LL(\Phi):=\lim_{n\to\infty}\frac1{n^2}\log E \exp(\Phi(\mmm))$$
is well defined for a large enough class of bounded continuous functions $\Phi:\calP\to\sR$,  and the proof of exponential tightness of $\{\mmm\}$.

\subsection{Varadhan functional}
Let $F_1,\dots F_m: \sR^2\to \sR$ be bounded continuous functions. Consider the bounded
continuous function $\Phi:\calP(\sR^2)\to\sR$ given by
\begin{equation}\lbl{Phi}
\Phi(\mu):=\min_{1\leq r\leq m} \int F_rd\mu.
\end{equation}
%Clearly, $\|\Phi\|_\infty\leq \max_r\|F_r\|_\infty<\infty$.

 We will show the following.

\begin{Theorem}\lbl{V-F}
Under the assumptions of Theorem \ref{T1},
\begin{eqnarray}\lbl{LL}
\lim_{n\to\infty} \frac1{n^2}\log \int \exp(n^2\Phi(\mmm))f(\xx,\yy)d\xx d\yy
\\ \nonumber
=\sup\left\{\Phi(\mu)-\int k(x,y)d\mu: \mu=\nu_1\otimes\nu_2\in \calP(\sR^2)\right\}.
\end{eqnarray}
\end{Theorem}
Denote $\KK(\mu)=\int k(x,y)d\mu$. Notice that by Assumption \ref{B1} we have
$\Phi(\mu)-\KK(\mu)\leq \max_r\|F_r\|_\infty+C$. In particular, 
\begin{equation}\lbl{finite}
\sup\{\Phi(\mu)-\KK(\mu): \mu\in \calP(\sR^2)\}<\infty.
\end{equation}
 
We prove \rf{LL} as two separate inequalities.  It will be convenient to prove the
upper bound for a larger class of functions $\Phi$. 
\begin{Lemma}\lbl{Upper bound}
If Assumptions \ref{p-int} and \ref{B1} hold true, then for every bounded continuous function
$\Phi:\calP(\sR^2)\to\sR$ we have
\begin{eqnarray}\lbl{LL-u}
\limsup_{n\to\infty} \frac1{n^2}\log \int \exp(n^2\Phi(\mmm))f(\xx,\yy)d\xx d\yy
\\ \nonumber
\leq\sup\left\{\Phi(\mu)-\KK(\mu): \mu=\nu_1\otimes\nu_2\in \calP(\sR^2)\right\}.
\end{eqnarray}
\end{Lemma}
\begin{proof} 
Notice that for $0<\theta<1$ 
$$\int \exp(n^2\Phi(\mmm))f(\xx,\yy)d\xx d\yy
$$
$$=\int \exp(n^2(\Phi(\mmm)-\theta \KK(\mmm))-(1-\theta)\sum_{i,j=1}^n k(x_i,y_j))d\xx d\yy
$$
$$
\leq\exp\left(n^2\sup_{\nu_1,\nu_2} (\Phi(\nu_1\otimes\nu_2)-
\theta \KK(\nu_1\otimes\nu_2))\right)\int \exp(-(1-\theta)\sum_{i,j=1}^n k(x_i,y_j))d\xx d\yy.
$$
Since 
$$\sum_{i,j=1}^n k(x_i,y_j)\geq -n^2C+\sum_{j=1}^nk(x_j,y_j),$$
 therefore
$$
\frac1{n^2}\log \int \exp(n^2\Phi(\mmm))f(\xx,\yy)d\xx d\yy$$
%$$
%\leq\sup_{\nu_1,\nu_2}\{\Phi(\nu_1\otimes \nu_2)-\theta K(\nu_1\otimes \nu_2)\}+(1-\theta)C+
%\frac1n\log\left(\iint g^{1-\theta}(x,y) dx dy\right)=
%$$
$$
\leq\sup_{\nu_1,\nu_2}\{\Phi(\nu_1\otimes \nu_2)-\theta K(\nu_1\otimes \nu_2)\}+(1-\theta)C+
\frac1n\log M_{1-\theta}.
$$
Thus
$$
\limsup_{n\to\infty} \frac1{n^2}\log \int \exp(n^2\Phi(\mmm))f(\xx,\yy)d\xx d\yy
$$
$$
\leq\sup_{\nu_1,\nu_2}\{\theta(\Phi(\nu_1\otimes \nu_2)-K(\nu_1\otimes
\nu_2))+(1-\theta)\Phi(\nu_1\otimes\nu_2)\}+2(1-\theta)C
$$
$$
\leq\theta \sup_{\nu_1,\nu_2}\{\Phi(\nu_1\otimes \nu_2)- K(\nu_1\otimes
\nu_2)\}+(1-\theta)\|\Phi\|_\infty+2(1-\theta)C.
$$
Passing to the limit as $\theta\to 1$ we get \rf{LL-u}.
\end{proof}

The proof of the lower bound is a combination of the discretization argument in 
\cite[pages 532--535]{BenArous-Guionnet-97} with the entropy estimate from  \cite[pages
191--192]{Johansson-98}. 
%It consists of a preliminary estimate, and an equivalent formula for the right-hand
%side of \rf{LL}.

Denote by $\calP_0$ the set of absolutely continuous probability measures $\nu(dx)=f(x)dx$
 on $\sR$ with compact support $\mbox{supp }(\nu)$,
and continuous density $f$.
Let us first record the well-known fact.
\begin{Lemma}\lbl{H}
If $\nu\in\calP_0$ then $\nu$ has finite entropy $$H_f:=\int \log f(x)\nu(dx)<\infty.$$   
\end{Lemma}
%\begin{proof} Suppose $\mbox{supp }( \nu) \subset [a,b]$. Since $f(x)$ is continuous, 
%$0\leq f(x)\leq C$ on $[a,b]$. Thus
%$$\inf_{0\leq u \leq C} u\log u \leq f(x)\log f(x) \leq \sup_{0\leq u \leq C} u\log u,$$ and
%$H_f=\int_a^b f(x) \log f(x) dx$ is finite. 
%\end{proof}

We first establish a weaker version of the lower bound.
\begin{Lemma}\lbl{Lower
bound}  If $\Phi$ is given by \rf{Phi}, then
\begin{eqnarray}\lbl{LL-low}
\liminf_{n\to\infty} \frac1{n^2}\log \int \exp(n^2\Phi(\mmm))f(\xx,\yy)d\xx d\yy
\\ \nonumber
\geq\sup\{\Phi(\nu_1\otimes \nu_2)-\KK(\nu_1\otimes \nu_2): \nu_1, \nu_2\in \calP_0\}.
\end{eqnarray}
\end{Lemma}
\begin{proof} 

Fix $\nu_1,\nu_2\in\calP_0$. Since $k(x,y)\geq -C$ is bounded from below, $\KK(\nu_1\otimes
\nu_2)\in(-\infty,\infty]$, so   without loss of generality we may assume that $k(x,y)$ is
$\nu_1\otimes \nu_2$-integrable.
 
 Since measures $\nu_1,\nu_2$ are absolutely continuous and have compact  supports, 
for every integer $n>0$ we can find  partitions 
$\Pi_1(n)=\{a_0< a_1<\dots< a_n\}$ and  $\Pi_2(n)=\{b_0< b_1<\dots< b_n\}$ of
$\mbox{supp }(\nu_1)$, $\mbox{supp }(\nu_2)$ respectively such that
$$\nu_1(a_{i-1},a_{i})=\nu_2(b_{j-1},b_{j})=\frac1n  \mbox{ for
$i,j=1,2,\dots,n.$} $$
Then, denoting 
$A=[a_0,a_1]\times[a_1,a_2]\times\dots\times[a_{n-1},a_n]$ and
$B=[b_0,b_1]\times[b_1,b_2]\times\dots\times[b_{n-1},b_n]$, we have 
\begin{eqnarray}\lbl{EF}
&\int \exp(n^2\Phi(\mmm))f(\xx,\yy)d\xx d\yy&
\\
\nonumber
&
%\int_{a_0}^{a_1}\int_{a_1}^{a_2}\dots \int_{a_{n-1}}^{a_n}
%\int_{b_0}^{b_1}\int_{b_1}^{b_2}\dots \int_{b_{n-1}}^{b_n}
\geq \int_{A\times B}
\exp\left(\min_r \sum_{i,j=1}^n F_r(x_i,y_j) - \sum_{i,j=1}^n k(x_i,y_j)\right) d\xx d\yy &
.\end{eqnarray}

Write $\nu_1=f(x)dx, \nu_2=g(y)dy$. By our choice of the partitions, functions
$$f_i(x):=nf(x)I_{[a_{i-1},a_i]}$$ and $$g_j(x):=ng(x)I_{[b_{j-1},b_j]}$$ are probability densities. 
Let 
$$
S(\xx,\yy)=\min_r \sum_{i,j=1}^n F_r(x_i,y_j) - \sum_{i,j=1}^n k(x_i,y_j)-\sum_{i=1}^n \log f(x_i)-
\sum_{j=1}^n \log g (y_j).
$$
Integrating over a smaller set $\{f_1(x_1)>0,\dots,f_n(x_n)>0,g_1(y_1)>0,\dots,g_n(y_n)>0\}$ on the right hand side of
\rf{EF}
we get $$
\int \exp(n^2\Phi(\mmm))f(\xx,\yy)d\xx d\yy$$
%$$
%\frac{1}{n^{2n}}\int
%\exp\left(\min_r \sum_{i,j} F_r(x_i,y_j) - \sum_{i,j} k(x_i,y_j)-\sum_{i} \log f(x_i)-\sum_{j} \log g (y_j)\right)$$
$$ 
\geq\frac{1}{n^{2n}}\int
\exp\left(S(\xx,\yy)\right)\prod_{i=1}^n f_i(x_i)\prod_{j=1}^n g_j(y_j)d\xx d\yy 
. $$
%(Here we omit a set of measure zero and integrate only over the set $\{f_1(x_1)>0,\dots,f_n(x_n)>0,g_1(y_1)>0,\dots,g_n(y_n)>0\}$).

Using Jensen's inequality, applied to the convex exponential function in the last integral, we get
$$
\int \exp(n^2\Phi(\mmm))f(\xx,\yy)d\xx d\yy\geq
\frac{1}{n^{2n}} \exp(S_1-S_2-S_3-S_4),
$$
where
\begin{eqnarray*}
S_1&=& \int_{\sR^{2n}} \left(\min_r \sum_{i,j=1}^n F_r(x_i,y_j)\right)\prod_{i=1}^n f_i(x_i)\prod_{j=1}^n g_j(x_j)d\xx
d\yy,
\\
S_2&=&\int_{\sR^{2n}}  \sum_{i,j=1}^n k(x_i,y_j)\prod_{i=1}^n f_i(x_i)\prod_{j=1}^n g_j(y_j)d\xx d\yy,
\\
S_3&=&
\int_{\sR^{n}} \sum_{i=1}^n \log f(x_i)\prod_{i=1}^n f_i(x_i)d\xx,
\\
S_4&=&\int_{\sR^{n}}  \sum_{j=1}^n \log g (y_j)\prod_{j=1}^n g_j(y_j) d\yy.
\end{eqnarray*}
 
We need the following identities. (Proofs of all Claims are postponed until the end of this proof.)
\begin{Claim}\lbl{int} For a  $\nu_1\otimes\nu_2$-integrable function $h$, we have
\begin{eqnarray*}
 \nonumber \int_{\sR^n} \sum_{j=1}^n h(y_j)\prod_{j=1}^n g_j(y_j) d\yy=n\int_{\sR}  h(y)g(y)dy,
\\
\int_{\sR^n} \sum_{i=1}^n h(x_i)\prod_{i=1}^n f_i(x_i) d\xx=n\int_{\sR}  h(x)f(x)dx,
\\
\int_{\sR^{2n}} \sum_{i,j=1}^n h(x_i,y_j)\prod_{i=1}^n f_i(x_i)\prod_{j=1}^n g_j(y_j) d\xx d\yy
\\
= n^2\iint  h(x,y)f(x)g(y)dxdy.
\end{eqnarray*}
\end{Claim}

Lemma \ref{H} says that the entropies $H_f=\int \log f(x) f(x) dx, H_g=\int \log g(y) g(y)
dy$ are finite. Thus the functions $k(x,y)$, $\log f(x)$, and  $\log g(y)$ are
$\nu_1\otimes\nu_2$-integrable.
Applying Claim \ref{int}, we get
$S_2= n^2 \KK(\nu_1\otimes\nu_2)$,
$S_3=n H_f$,
and
$S_4=
n H_g$. Therefore,
\begin{eqnarray}\lbl{*}
&\int \exp(n^2\Phi(\mmm))f(\xx,\yy)d\xx d\yy &\\ \nonumber &\geq 
\frac{1}{n^{2n}} \exp(S_1-n^2 \KK(\nu_1\otimes\nu_2) -n H_f-n H_g)&.
\end{eqnarray}

We  need the following lower bound for $S_1$.
\begin{Claim}\lbl{min}
\begin{eqnarray}\lbl{mmm}
\int \left(\min_r \sum_{i,j=1}^n F_r(x_i,y_j)\right)\prod f_i(x_i)\prod g_j(y_j)d\xx d\yy
\\ \nonumber
\geq\min_r \sum_{i,j=1}^n F_{r,(i,j)},
\end{eqnarray}
where 
$$F_{r,(i,j)}=\min \left\{F_r(x,y): a_{i-1}\leq x\leq a_i, b_{j-1}\leq y\leq b_j\right\}.$$

\end{Claim}
 
Combining inequalities \rf{*} and \rf{mmm}, we get
\begin{eqnarray}\lbl{LU}
\frac1{n^2}\log \int \exp(n^2\Phi(\mmm))f(\xx,\yy)d\xx d\yy
\\ \nonumber
\geq\frac1{n^2}\min_r \sum_{i,j=1}^n F_{r,(i,j)}- \KK(\nu_1\otimes\nu_2)
-\frac1n H_f - \frac1n H_g -\frac2n\log n.
\end{eqnarray}
%where $H_f=\int \log f(x) f(x) dx, H_g=\int \log g(y) g(y) dx$ are finite.

Since functions $F_r(x,y)$ are continuous and $\nu_1,\nu_2$ have compact support and continuous densities $f$, $g$,  therefore  
$\nu_1\otimes \nu_2$-almost surely 
$\sum_{i,j=1}^n F_{r,(i,j)}I_{(a_{i-1},a_i)}(x)I_{(b_{j-1},b_j)}(y)\to F_r(x,y)$ (to see this, notice that  for fixed $\eps>0$, the sequence convergences for all $(x,y)$ such that 
$f(x)\geq \eps, g(y)\geq \eps$.),  
and the functions are bounded. 
Since $1\leq r\leq m$ ranges over a finite set of
values only we have
$$\lim_{n\to\infty} \frac1{n^2}\min_r \sum_{i,j=1}^n F_{r,(i,j)}$$ $$
=\min_r \lim_{n\to\infty}\sum_{i,j=1}^n F_{r,(i,j)}\nu_1(a_{i-1},a_i)\nu_2(b_{j-1}b_j)
=\Phi(\nu_1\otimes\nu_2).$$
Letting  $n\to\infty$ in \rf{LU}  we obtain \rf{LL-low}.

To conclude the proof, it remains to prove Claims \ref{int} % 
and \ref{min}.
 
\begin{proofof}{Proof of Claim {\protect\ref{int}}} 
Switching the order of integration and summation, we get 
$$\int \sum_{j=1}^n h(y_j)\prod g_j(y_j) d\yy
=\sum_{j=1}^n \int h(y_j)g_j(y_j)dy_j\prod_{i\ne j}\int g_i(y_i)dy_i$$
$$
=\sum_{j=1}^n \int h(y_j)g_j(y_j)dy_j=
n\sum_{j=1}^n \int_{b_{j-1}}^{b_j} h(y)g(y)dy
=
n \int_{b_{0}}^{b_n} h(y)g(y)dy.
$$
The other two identities follow by a similar argument.
% for the first one we switch the roles
%of $\xx,\yy$; the last identity follows from the first two: consider    $h_1(y):=\sum_i
%h(x_i,y)$ as a function of $y$ with fixed $\xx$, and then integrate with respect to $\xx$.
\end{proofof}

\begin{proofof}{Proof of Claim {\protect\ref{min}}} 
Fix $0\leq  k\leq n$, $x_1,\dots,x_{k}\in\sR$ and $y_1,\dots,y_n\in\sR$. Let 
 $$G_{r,k}(x_1,\dots,x_{k}):=\sum_{i=1}^k \sum_{j=1}^n F_r(x_i,y_j)+\sum_{i=k+1}^n
\sum_{j=1}^n\min_{a_{i-1}\leq x\leq a_i}F_r(x,y_j).$$

If $a_{k-1}<x_k<a_k$, we have
 $$\min_rG_{r,k}(x_1,\dots,x_{k})=
$$
$$
\min_r\left(\sum_{i=1}^{k-1} \sum_j F_r(x_i,y_j)+\sum_j F_r(x_k,y_j)+
\sum_{i=k+1}^n \sum_j\min_{a_{i-1}\leq x\leq a_i}F_r(x,y_j)\right)
$$
$$
\geq\min_r G_{r,k-1}(x_1,\dots,x_{k-1}).
$$
Therefore,
$$
\int_{a_{k-1}}^{a_{k}}
\min_r G_{r,k}(\xx)f_{k}(x_{k}) d x_k 
\geq
\min_r G_{r,k-1}(\xx).
$$
Recurrently, 
$$
\int \min_r \left(\sum_{i=1}^{n}\sum_{j=1}^{n} F_r(x_i,y_j)\right)\prod f_i(x_i)d\xx 
$$
$$
=\int \min_r G_{r,n}(\xx)\prod f_i(x_i)d\xx \geq \min_r G_{r,0}(\xx)
$$
$$
=\min_r \left(\sum_{i=1}^{n}\sum_{j=1}^{n} \min_{a_{i-1}\leq x\leq a_i}F_r(x,y_j)\right).
$$
Applying the same reasoning to variables $y_1,\dots,y_n$ and 
$$G_{r,k}(y_1,\dots,y_k):=
\sum_{j=1}^{k}\sum_{i=1}^{n}\min_{a_{i-1}\leq x\leq a_i}F_r(x,y_j)+\sum_{j=k+1}^n
\sum_{i=1}^{n}F_{r,(i,j)}
$$ 
we get \rf{mmm}.
\end{proofof}

This concludes the proof.
\end{proof}

%\subsubsection{Reductions}
The next Lemmas  show that the right hand sides of \rf{LL-u} and \rf{LL-low}  coincide.

 Let $\calP_c$ denote 
compactly supported probability measures.
\begin{Lemma}\lbl{comp} If Assumption \ref{B1} holds true, then
\begin{eqnarray}\lbl{comp *}
\sup\{\Phi(\nu_1\otimes \nu_2)-\KK(\nu_1\otimes \nu_2): \nu_1, \nu_2\in \calP_c\}
\\ \nonumber 
=\sup\{\Phi(\nu_1\otimes \nu_2)-\KK(\nu_1\otimes \nu_2): \nu_1, \nu_2\in \calP\}.
\end{eqnarray}
\end{Lemma}
\begin{proof} Clearly the left-hand side of \rf{comp *} cannot exceed the right hand side. 
To show the converse inequality, fix $\eta>0$ and  $\nu_1, \nu_2\in \calP$
such that 
\begin{equation}\lbl{*0}
\Phi(\nu_1\otimes \nu_2)-\KK(\nu_1\otimes \nu_2)\geq  \sup_{\nu_1, \nu_2\in
\calP}\{\Phi(\nu_1\otimes
\nu_2)-\KK(\nu_1\otimes \nu_2) \} -\eta.
\end{equation}
Since the supremum is finite, see \rf{finite}, and $k(x,y)$ is bounded from below,
therefore $\iint |k(x,y)| d\nu_1d\nu_2<\infty$.
 
For $L>0$ large enough, define probability measures 
$\nu_{j,L}$ by
$$\nu_{j,L}(A):=\frac{\nu_j(A\cap[-L,L])}{\nu_j([-L,L])}, \, j=1,2.$$
By definition, measures $\nu_{j,L}\in\calP_c$ have compact support.
Since $-C\leq k(x,y)I_{|x|<L, |y|<L}\leq |k(x,y)|$ and $k$ is
$\nu_1\otimes\nu_2$-integrable, by  Lebesgue's dominated convergence theorem
$$\lim_{L\to\infty}\KK(\nu_{1,L}\otimes\nu_{2,L})=
\frac{\lim_{L\to\infty}\int_{-L}^L\int_{-L}^L k(x,y)\nu_1(dx)\nu_2(dy)}
{\lim_{L\to\infty}\nu_1([-L,L])\nu_2([-L,L])}=
\KK(\nu_1\otimes \nu_2).$$
 Similarly, 
$$\lim_{L\to\infty}\Phi(\nu_{1,L}\otimes\nu_{2,L})=\Phi(\nu_1\otimes \nu_2).$$
Thus \rf{comp *} follows. 
\end{proof}

\begin{Lemma} If Assumptions \ref{B1} and \ref{B3} hold true, then
\begin{eqnarray}\lbl{reduct}
\sup\{\Phi(\nu_1\otimes \nu_2)-\KK(\nu_1\otimes \nu_2): \nu_1, \nu_2\in \calP_0\}
\\ \nonumber 
=\sup\{\Phi(\nu_1\otimes \nu_2)-\KK(\nu_1\otimes \nu_2): \nu_1, \nu_2\in \calP\}.
\end{eqnarray}
\end{Lemma}
\begin{proof} Trivially,
$$
\sup_{\nu_1, \nu_2\in \calP_0}\{\Phi(\nu_1\otimes \nu_2)-\KK(\nu_1\otimes \nu_2) \}
\leq
\sup_{\nu_1, \nu_2\in \calP}\{\Phi(\nu_1\otimes \nu_2)-\KK(\nu_1\otimes \nu_2) \}.$$
To show the converse inequality, fix $\eta>0$ and compactly supported  $\nu_1, \nu_2\in
\calP_c$ such that 
\begin{equation}\lbl{*1}
\Phi(\nu_1\otimes \nu_2)-\KK(\nu_1\otimes \nu_2)\geq  \sup_{\nu_1, \nu_2\in
\calP}\{\Phi(\nu_1\otimes
\nu_2)-\KK(\nu_1\otimes \nu_2) \} -\eta,
\end{equation}
see Lemma \ref{comp}. As previously,  $k(x,y)$ is
$\nu_1\otimes\nu_2$-integrable, see \rf{finite}.

Consider the convolution $\nu_{j}^{\eps}(A):=\frac{1}{2\eps}\int_{-\eps}^\eps\nu_j(A-x)dx$,
where $j=1,2$ and  $0<\eps\leq 1$.  Measures $\nu_{1}^{\eps},\nu_{2}^{\eps}$ have continuous
densities, and  since $\nu_1,\nu_2$ have compact supports, $\nu_{1}^{\eps},\nu_{2}^{\eps}$
also have compact support. Thus $\nu_{1}^{\eps},\nu_{2}^{\eps}\in\calP_0$ and
\begin{equation}\lbl{*00}
\sup_{\nu_1, \nu_2\in \calP_0}\{\Phi(\nu_1\otimes \nu_2)-\KK(\nu_1\otimes \nu_2) \}\geq 
\Phi(\nu_{1}^{\eps}\otimes\nu_{2}^{\eps})-\KK(\nu_{1}^{\eps}\otimes\nu_{2}^{\eps}).
\end{equation}

 As $\eps\to 0$ measure $\nu_{j}^{\eps}$ converges weakly to $\nu_{j}$. Hence
\begin{equation}\lbl{**}
\lim_{\eps\to0}\Phi(\nu_{1}^{\eps}\otimes\nu_{2}^{\eps})=\Phi(\nu_1\otimes\nu_2).
\end{equation}

Assumption \ref{B3} asserts that $V(x,y):=\beta\log|x-y|+k(x,y)$ is a continuous function.  
Thus $|V(x,y)|$ is bounded on the compact set $\mbox{supp }(\nu_1^1\otimes\nu_2^1)$.
Since the supports of  $\nu_{1}^{\eps}\otimes\nu_{2}^{\eps}$ are 
contained in $\mbox{supp }(\nu_1^1\otimes\nu_2^1)$, and $\nu_{j}^{\eps}\to\nu_{j}$, we get
\begin{equation}\lbl{***}
\iint V(x,y)\nu_1^\eps(dx)\nu_2^\eps(dy)\to \iint V(x,y)\nu_1(dx)\nu_2(dy).
\end{equation}
This concludes the proof if $\beta=0$. If $\beta>0$, then $\log|x-y|$ is
$\nu_1\otimes\nu_2$-integrable as a linear combination of
integrable functions, $\log|x-y|=(V(x,y) -k(x,y))/\beta$. Therefore we have
$$
\Phi(\nu_{1}^{\eps}\otimes\nu_{2}^{\eps})-\KK(\nu_{1}^{\eps}\otimes\nu_{2}^{\eps})=
\Phi(\nu_{1}^{\eps}\otimes\nu_{2}^{\eps})$$
$$
-\iint V(x,y)\nu_1^\eps(dx)\nu_2^\eps(dy)+
\beta\iint\log|x-y|\nu_1(dx)\nu_2(dy)$$
$$-\beta\left(\iint \log|x-y|\nu_1(dx)\nu_2(dy)-\iint \log|x-y|\nu_1^\eps(dx)\nu_2^\eps(dy)\right).
$$ 
Taking the lim sup as $\eps\to 0$, from \rf{*00}, \rf{**}, \rf{***},  and  \rf{*1} we get
$$
\sup_{\nu_1, \nu_2\in \calP_0}\{\Phi(\nu_1\otimes \nu_2)-\KK(\nu_1\otimes \nu_2) \}\geq
\sup_{\nu_1, \nu_2\in \calP}\{\Phi(\nu_1\otimes
\nu_2)-\KK(\nu_1\otimes \nu_2) \} -\eta
$$
$$
 -\limsup_{\eps\to 0}  
\left(\iint \log|x-y|\nu_1(dx)\nu_2(dy)-\iint \log|x-y|\nu_1^\eps(dx)\nu_2^\eps(dy)\right).
$$
Since $\eta>0$ is arbitrary, to end the proof   we  use  the following. 
\begin{Claim}\lbl{log-est}
$$\limsup_{\eps\to 0} 
\left(\iint \log|x-y|\nu_1(dx)\nu_2(dy)-\iint \log|x-y|\nu_1^\eps(dx)\nu_2^\eps(dy)\right)\leq 0.$$
\end{Claim}
\end{proof}

\begin{proofof}{Proof of Claim {\protect\ref{log-est}}}
Claim \ref{log-est} is established by the argument in \cite[pages 192-193]{Johansson-98}. For completeness,
 we repeat it here. %\mybox{Really? Or did he have one measure only?}
Let $X,Y$ be independent random variables with distributions $\nu_1,\nu_2$ respectively and
let $Z=X-Y$. Since $\log|Z|$ is integrable, $\Pr(Z=0)=0$.  Let $U\in[-2,2]$
be a r.~v. independent of $Z$  with the density $f(u)=(2-|u|)/4$.  It is easy to see that
$\iint \log|x-y|\nu_1^\eps(dx)\nu_2^\eps(dy)=E\log |Z+\eps U|$, and the inequality to prove
reads
$$
\limsup_{\eps\to 0}E\left(\log^+\frac{1}{|1+\eps \frac{U}{Z}|}\right)\leq 0.
$$
For fixed $z\ne 0$ we have 
\begin{equation}\lbl{Joh}
E\left(\log^+\frac{1}{|1+\eps \frac{U}{z}|}\right)\leq \frac1{\log
2}\log\left(1+\frac{2\eps}{|z|}\right).
\end{equation}
Indeed, since $(2-|u|)/4\leq 1/2$ we get
$$
E\left(\log^+\frac{1}{|1+\eps \frac{U}{z}|}\right)\leq 
\frac{|z|}{4\eps}\int_{1-2\eps/|z|}^{1+2\eps/|z|}\log^+\frac{1}{|x|} dx.$$
Therefore,
$$
E\left(\log^+\frac{1}{|1+\eps \frac{U}{z}|}\right)\leq \left\{\begin{array}{ll}
\frac{|z|}{4\eps}\int_{1-2\eps/|z|}^{1}\log\frac{1}{x} dx & \mbox{ if $|z|>2\eps$}\\
\frac{|z|}{4\eps}\left(\int_0^1\log \frac{1}{x} dx +\int_{0}^{2\eps/|z|-1}\log^+\frac{1}{x} dx\right) 
& \mbox{ if $|z|\leq 2\eps$}
\end{array}
\right. .
$$
If $|z|>2\eps$ we get 
$E\left(\log^+\frac{1}{|1+\eps \frac{U}{z}|}\right)\leq 
%\frac{|z|}{4\eps} \frac{2\eps}{|z|}\log \frac{1}{1-2\eps/|z|}=
\frac{1}{2} \log \frac{1}{1-2\eps/|z|} 
%=\frac12\log\left(1+\left(\frac{2\eps}{|z|}\right)+\left(\frac{2\eps}{|z|}\right)^2+\dots\right)
<\frac12\log(1+\frac{2\eps}{|z|})\leq \frac1{\log 2}\log(1+\frac{2\eps}{|z|})
$.
If $|z|\leq 2\eps$, then
$
E\left(\log^+\frac{1}{|1+\eps \frac{U}{z}|}\right)\leq 
\frac{|z|}{2\eps}\int_0^1\log \frac{1}{x} dx\leq 1 \leq 
\frac1{\log 2}\log(1+\frac{2\eps}{|z|})
$.  Thus in both cases, \rf{Joh} holds true.

To finish the proof we integrate inequality \rf{Joh} and get
$$
\limsup_{\eps\to 0}E\left(\log^+\frac{1}{|1+\eps \frac{U}{Z}|}\right)\leq
\frac1{\log 2}\limsup_{\eps\to 0}E\left(\log(1+2\eps /|Z|)\right).
$$ 
For $\eps<1/2$ we have 
$\log(1+2\eps /|Z|)\leq \log 2+ \log^+\frac1{|Z|}$ and $\log^+\frac1{|Z|}$ is integrable.
Lebesgue's dominated convergence theorem yields
$$\limsup_{\eps\to 0}E\left(\log(1+2\eps /|Z|)\right)=0.$$
\end{proofof} 

\begin{proofof}{Proof of Theorem \protect{\ref{V-F}}}
Combining Lemmas \ref{Upper bound} and \ref{Lower bound} we have
$$
\sup\{\Phi(\nu_1\otimes \nu_2)-\KK(\nu_1\otimes \nu_2): \nu_1, \nu_2\in \calP_0\}
$$
$$
\leq\liminf_{n\to\infty} \frac1{n^2}\log \int \exp(n^2\Phi(\mmm))f(\xx,\yy)d\xx d\yy
$$
$$
\leq\limsup_{n\to\infty} \frac1{n^2}\log \int \exp(n^2\Phi(\mmm))f(\xx,\yy)d\xx d\yy
$$ 
$$
\leq\sup\{\Phi(\nu_1\otimes \nu_2)-\KK(\nu_1\otimes \nu_2): \nu_1, \nu_2\in \calP\}.
$$
By \rf{reduct}, all of the above inequalities are in fact equalities. Thus \rf{LL} holds true. 
\end{proofof}
\subsection{Exponential tightness}
Recall that $\{\mmm\}$ is exponentially tight if
for every $m>0$ there is a compact subset $\calK\subset\calP$ such that
$$
\sup_n\frac1{n^2}\log \Pr(\mmm\not\in \calK)<-m.
$$
Our proof of exponential tightness is a concrete implementation of de Acosta \cite{deAcosta-85}.
%, see also \cite{Deuschel-Stroock}, \cite{Dembo-Zeitouni}, \cite{Dupuis-Ellis}.

Assumption \ref{B1} implies that $k(x,y)+C\geq 0$.
Let $q:\calP(\sR^2)\to [0,\infty]$ be given by $$q(\mu)=\int_{\sR^2} \left(k(x,y)+C\right)d\mu.$$ 
\begin{Lemma}\lbl{pre-comp} If Assumptions \ref{B1}, \ref{B2} hold true, then
$q$ has pre-compact level sets:  for every $t>0$, $q^{-1}[0,t]$ is a pre-compact set in $\calP$.
\end{Lemma}  
\begin{proof} Fix $t>0$ and denote $\calK:=\{\mu\in\calP(\sR^2): q(\mu)\leq t\}$. We will
show that $\calK$ is pre-compact. 

Assumption \ref{B2} says that for every $\eps>0$ the set
$K_\eps:=\{(x,y): C+k(x,y)\leq t/\eps\}$ is a compact subset of $\sR^2$. For every $\mu\in \calK$
 by Chebyshev's
inequality we have $$\mu(K_\eps^c)\leq \mu(\{(x,y):C+k(x,y)>t/\eps\})\leq \frac{\eps q(\mu)}{t} =\eps.$$ 
Thus $\calK$ is pre-compact,
and its weak closure $\bar{\calK}$ is compact.
\end{proof}
 \begin{Lemma}\lbl{exp(q)} 
If Assumptions \ref{B0}, \ref{B1}, and \ref{p-int}  hold true, then 
\begin{equation}\lbl{q2}
\sup_n \frac1{n^2}\log \int\exp(\frac12n^2 q(\mmm))f(\xx,\yy)d\xx d\yy <\infty.
\end{equation}
\end{Lemma}
\begin{proof} 
We have $$\int\exp(\frac12n^2 q(\mmm))f(\xx,\yy)d\xx d\yy=
\int \exp( \frac12n^2C -\frac12\sum_{i,j=1}^nk(x_i,y_j) )d\xx d\yy
 $$
$$
\leq e^{\frac12n^2C}\int \prod_{i,j=1}^n \sqrt{g(x_i,y_j)} d\xx d\yy 
\leq
%e^{\frac12n^2C}\left(\int\sqrt{g(x,y)e^{C(n-1)}} dx dy\right)^n \leq
e^{n^2C}M_{1/2}^{n}.
$$
Therefore the left-hand side of \rf{q2} is at most $C+\log^+M_{1/2} <\infty$.
\end{proof}
\begin{Theorem}\lbl{Exp-T}
Under the assumptions of Theorem \ref{T1}, the sequence $\{\mmm\}$ is exponentially tight.
\end{Theorem}
\begin{proof}

Notice that by \rf{LL} used with $\Phi(\mu):=0$ we have
$\frac{1}{n^2}\log Z_n\to L_0:=-\inf_\mu \int k(x,y) d\mu =-\inf_{x,y}k(x,y)$.
Since $L_0$ is finite, see \rf{finite}, therefore by Lemma \ref{exp(q)} we have
$$
\sup_n \frac1{n^2}\log \int\exp(\frac12n^2 q(\mmm))\frac{1}{Z_n}f(\xx,\yy)d\xx d\yy =C_1<\infty.
$$
Fix $m>0$. Let $\calK\subset \calP$ be the pre-compact set from Lemma \ref{pre-comp}
corresponding to $t=2m+2C_1$.

Applying Chebyshev's inequality to probability measure \rf{Pr}
% $$\Pr(d\xx,d\yy):=\frac1{Z_n}f(\xx,\yy)d\xx d\yy$$ 
we get
$$\Pr(\mmm\not\in\bar{\calK})\leq \Pr(\mmm\not\in \calK)=\Pr(q(\mmm)>t)
\leq e^{-\frac12n^2t} \int \exp(\frac12n^2q(\mmm))d\Pr.
$$
Therefore $$\Pr(\mmm\not\in\bar{\calK})\leq  e^{-\frac12n^2t} e^{n^2 C_1},$$ and
$$\frac1{n^2}\log \Pr(\mmm\not\in\bar{\calK})\leq -t/2+C_1=-m 
$$
for all $n$.
\end{proof}
%\subsection{Proof of Theorem {\protect\ref{T1}}}
\begin{proofof}{Proof of Theorem {\protect\ref{T1}}}
Recall that the space $\calP=\calP(\sR^2)$ of probability measures on $\sR^2$ with the
topology of weak convergence is a Polish space.  By Theorem \ref{Exp-T}, $\{\mmm\}$ is
exponentially tight. Theorem \ref{V-F} says that the Varadhan functional 
$$\LL(\Phi):=\lim_{n\to\infty} \frac1{n^2}\log E \left(\exp n^2\Phi(\mmm)\right)$$
is defined on all functions $\Phi$ given by \rf{Phi}, and
$$\LL(\Phi)=\sup_{\nu_1,\nu_2}\{\Phi(\nu_1\otimes \nu_2)-\KK(\nu_1\otimes
\nu_2)\}-\lim_{n\to\infty}\frac1{n^2}\log Z_n  $$
$$
=\sup_{\nu_1,\nu_2}\{\Phi(\nu_1\otimes \nu_2)-\KK(\nu_1\otimes
\nu_2)\}+\inf_{x,y}k(x,y).$$
Thus
\begin{equation}\lbl{@}\LL(\Phi)=\sup_{\nu_1,\nu_2}\{\Phi(\nu_1\otimes \nu_2)-\KK(\nu_1\otimes
\nu_2)\}+I_0.\end{equation}

Functions $\Phi$ defined by \rf{Phi} form a subset of 
$C_b(\calP(\sR^2))$  which separates points of $\calP(\sR^2)$ and is closed under the operation of taking
pointwise
minima.
 Thus by \cite[Theorem T.1.3]{Bryc90a} or \cite[Theorem 4.4.10]{Dembo-Zeitouni},
the empirical measures $\{\mmm\}$ satisfy the large deviation principle with the rate function
\begin{equation}\lbl{III}
\II(\mu):=\sup\{\Phi(\mu)-\LL(\Phi)\};
\end{equation}
here, the supremum is taken over all  $F_1,\dots,F_m\in C_b(\sR^2)$ and $\Phi(\mu)$ is defined by 
\rf{Phi}. 

It remains to prove formula \rf{I}. Fix $\nu_1,\nu_2\in\calP$.
From \rf{@}, for $\Phi$ given by \rf{Phi} we have
$\LL(\Phi)%=\sup_{\nu_1, \nu_2\in \calP}\{\Phi(\nu_1\otimes\nu_2)-\KK(\nu_1\otimes \nu_2) \}+I_0 
\geq \Phi(\nu_1\otimes \nu_2)-\KK(\nu_1\otimes \nu_2)+I_0$. 
Thus formula \rf{III} implies that
\begin{equation}\lbl{I+}
\II(\nu_1\otimes \nu_2)\leq \KK(\nu_1\otimes
\nu_2)-I_0.
\end{equation}

To prove the converse inequality we use the fact that we already know that 
the large deviations principle holds.
The large deviations principle implies that
\begin{equation}\lbl{LDP}
\II(\nu_1\otimes\nu_2)=\sup_{\Phi\in C_b(\calP)}\{\Phi(\nu_1\otimes\nu_2) -\LL(\Phi)\}.
\end{equation} 

Now consider $\Phi_M(\mu)=\int (M\wedge k(x,y))\, d\mu$. Assumptions \ref{B1} and \ref{B3}
imply that $(x,y)\mapsto M\wedge k(x,y)$ is a bounded continuous function for every real
$M$. %(Indeed, $k(x,y)=V(x,y)-\beta\log|x-y|$ for a continuous function $V$.) 
Thus
$\Phi_M$ is given by \rf{Phi}.
 Since $M\wedge k(x,y)\leq k(x,y)$, from %\rf{LL} and 
\rf{@} we get
$\LL(\Phi_M)\leq I_0$. 
Thus 
$$\II(\nu_1\otimes\nu_2)\geq \sup_{M}\{\Phi_M(\nu_1\otimes\nu_2) -\LL(\Phi_M)\}\geq
\limsup_{M\to\infty}\int M\wedge k(x,y)d\mu -I_0.$$
This together with \rf{I+} proves \rf{I} for product measures.

%
%Inded, $\II(\nu_1\otimes \nu_2)\geq  \Phi(\nu_1\otimes
%\nu_2)-\LL(\Phi)$ which is the same as 
% $\LL(\Phi)=\sup_\mu ... \geq \Phi(\nu_1\otimes
%\nu_2) -\II(\nu_1\otimes \nu_2)$. On the otehr hand 
%$\II(\nu_1\otimes \nu_2)\leq  \Phi(\nu_1\otimes
%\nu_2)-\KK(\nu_1\otimes \nu_2) \}-I_0$. 

It remains to verify that if $\mu_0$ is not a product measure, then 
$\II(\mu_0)=\infty$.  To this end, take  bounded continuous functions  $F(x),G(y)$
such that  $$\delta:=\int F(x)G(y)\mu_0(dx,dy)- \int F(x)\mu_0(dx,dy)\int G(y)\mu_0(dx,dy)>0.$$
For $b>0$, let
$$\Phi_b(\mu):=b\left(\int F(x)G(y)\mu(dx,dy)- \int F(x)\mu(dx,dy)\int G(y)\mu(dx,dy)\right).
$$ 
 
Clearly, $\Phi_b:\calP\to\sR$ is a bounded continuous function, which vanishes on
product measures. By the upper bound \rf{LL-u} we therefore have
$\LL(\Phi_b)\leq I_0$. So
$\II(\mu_0)\geq \Phi_b(\mu_0)-\LL(\Phi_b)\geq b\delta-I_0$.
Since $b$ can be arbitrarily large, $\II(\mu_0)=\infty$.
\end{proofof}

\subsection*{Acknowledgements} 
I would like to thank P. Dupuis for a conversation on non-convex rate functions.

\def\cprime{$'$} \def\cprime{$'$} \def\cprime{$'$} \def\cprime{$'$}
  \def\cprime{$'$}

\end{document}